\documentclass[12pt]{amsart}
\newtheorem{theorem}{Theorem}
\newtheorem{proposition}[theorem]{Proposition}
\def\ds{\displaystyle}

\title[]
{The symmetrized polydisc cannot be exhausted by domains
biholomorphic to convex domains}

\author{Nikolai Nikolov}

\address
{Institute of Mathematics and Informatics\\ Bulgarian Academy of
Sciences\\ 1113 Sofia, Bulgaria} \email{nik@math.bas.bg}

\subjclass[2000]{32H35}

\keywords{symmetrized polydisc}

\begin{document}

\begin{thanks}
{This paper was started during the stay of the author at the
Jagiellonian University, Krak\'ow in February, 2005. He would like
to thank Peter Pflug, Marek Jarnicki and W\l odzimierz Zwonek for
helpful discussions.}
\end{thanks}

\begin{abstract} We prove that the symmetrized polydisc cannot be exhausted
by domains biholomorphic to convex domains.
\end{abstract}

\maketitle

Let $\Bbb D$ be the unit disc in $\Bbb C.$ Let
$\sigma_n=(\sigma_{n,1},\dots,\sigma_{n,n}):\Bbb C^n\to\Bbb C^n$
be defined as follows: $$\sigma_{n,k}(z_1,\dots,z_n)=\sum_{1\le
j_1<\dots<j_k\le n}z_{j_1}\dots z_{j_k},\ \ 1\le k\le n.$$ The set
$\Bbb G_n=\sigma_n(\Bbb D^n)$ is called the symmetrized $n$-disc.
The symmetrized bidisc $\Bbb G_2$ is the first example of a
bounded pseudoconvex domain, which is not biholomorphic to any
convex domain and on which the Carath\'eodory and Kobayashi
distances coincide (see \cite{Cos}). Moreover, it cannot be
exhausted by domains biholomorphic to convex domains (see
\cite{Edi}). It has been asked in \cite{Jar-Pfl2} whether the last
result remains true for $\Bbb G_n,$ $n\ge 3.$ The aim of this note
is to give a positive answer to the above question.

Let us begin with the following definition. Let $k_1\le\dots\le
k_n$ be positive integers and
$$\pi_\lambda(z_1,\dots,z_n)=(\lambda^{k_1}z_1,\dots,\lambda^{k_n}z_n).$$
A domain $D$ in $\Bbb C^n$ is called $(k_1,\dots,k_n)$-balanced if
$\pi_\lambda(z)\in D$ for $z\in D,\lambda\in\overline{\Bbb D}$.
For such a domain $D$ one has $$D=\{z\in\Bbb C^n:h(z)<1\},$$ where
$$h(z)=\inf\{(\lambda>0:\pi_{1/\lambda}(z)\in D\},\ z\in\Bbb
C^n.$$ It is easy to see that $h$ is an upper semicontinuous,
non-negative function on $\Bbb C^n$ with
$$h(\pi_\lambda(z))=|\lambda|h(z),\ \lambda\in\Bbb C,z\in\Bbb
C^n.$$

Note that the $(1,\dots,1)$-balanced domains are exactly the
balanced domains in the usual sense (cf. \cite{Jar-Pfl1}). As in
the case of balanced domains one has the following

\begin{proposition} A $(k_1,\dots,k_n)$-balanced domain $D$ is pseudoconvex
if and only $\log h$ is a plurisubharmonic function.
\end{proposition}

\begin{proof} It is clear that if $\log h$ is a plurisubharmonic
function, then $D$ is a pseudoconvex domain.

To prove the converse, define $\Phi:\Bbb
C^n\owns(z_1,\dots,z_n)\mapsto (z_1^{k_1},\ldots,z_n^{k_n})\in\Bbb
C^n$ and set $\tilde D:=\Phi^{-1}(D)$, $\tilde h=h\circ\Phi$. Note
that $\tilde D=\{z\in\Bbb C^n:\tilde h(z)<1\}$ and $\tilde
h(\lambda z)=|\lambda|h(z)$, $\lambda\in\Bbb C$, $z\in\Bbb C^n$.
Therefore $\tilde D$ is a pseudoconvex balanced domain whose
Minkowski functional is equal to $\tilde h$. Consequently,
$\log\tilde h$ is a plurisubharmonic function (cf.
\cite{Jar-Pfl1}). On the other hand, one has $h(z)=\tilde
h(\root{k_1}\of{z_1},\ldots,\root{k_n}\of{z_n})$, $z\in\Bbb
C_*^n$, where the roots are arbitrarily chosen. Thus $\log h$ is a
plurisubharmonic function on $\Bbb C_*^n$ and hence, by the
removable singularities theorem (cf. \cite{Jar-Pfl1}), it is
plurisubharmonic on $\Bbb C^n$.
\end{proof}

The crucial step in the proof of our main result is the following

\begin{proposition} Let $D$ be a $(k_1,\dots,k_n)$-balanced
domain, which can be exhausted by domains biholomorphic to convex
domains. If $2k_{m+1}>k_n$ for some $m,$ $0\le m\le n-1,$ then the
intersection $D_m=D\cap\{z_1=\dots=z_m=0\}$ is a convex set (we
assume that $D_m=D$ if $m=0$).
\end{proposition}

\begin{proof} The proof is similar to that of Theorem 1 in
\cite{Edi}.

Take two points $a,b\in D_m.$ We may find a domain $D'\subset D$
which is biholomorphic to a convex domain $G$ and such that
$\lambda a,\lambda b\in D$ for $\lambda\in\overline{\Bbb D}.$ Let
$\Psi:D'\to G$ be the corresponding biholomorphic mapping. We may
assume that $\Psi(0)=0$ and $\Psi'(0)=\hbox{id}.$ If
$$g_{ab}(\lambda)=\frac{\Psi(\pi_\lambda(a))+\Psi(\pi_\lambda(b))}{2},$$
then $\Psi^{-1}\circ g_{ab}(\lambda)$ is a holomorphic mapping
from a neighborhood of $\overline{\Bbb D}$ into $D.$ Set
$f_{ab}(\lambda)=\pi_{1/\lambda}\circ\Psi^{-1}\circ
g_{ab}(\lambda).$ We shall see later that $f_{ab}(\lambda)$ can be
extended at $0$ by proving that $$\lim_{\lambda\to
0}f_{ab}(\lambda)=\frac{a+b}{2}.\eqno{(1)}$$ If (1) holds, then
$h\circ f_{ab}$ is a subharmonic function by Proposition 1 and the
maximum principle implies that
$$h(f_{ab}(0))\le\max_{|\lambda|=1}h(f_{ab}(\lambda))<1.$$ Hence
$\ds\frac{a+b}{2}\in D_m$ if $a,b\in D_m,$ i.e. $D_m$ is a convex
set.

To prove (1), note that $\Psi^{-1}(0)=0$ and
$(\Psi^{-1})'(0)=\hbox{id}$ imply that, for any $j=1,2,\dots,n,$
one has $$\Psi_j^{-1}\circ
g_{ab}(\lambda)=g_{abj}(\lambda)+O(|g_{ab}(\lambda)|^2).$$ Since
$\Psi(0)=0,$ $\Psi'(0)=\hbox{id}$ and $a,b\in D_m,$ it follows
that
$$g_{abj}(\lambda)=\frac{a_j+b_j}{2}\lambda^j+O(|\lambda|^{2k_{m+1}}).$$
Now the inequality $2k_{m+1}>k_n$ shows that
$$\frac{\Psi_j^{-1}\circ
g_{ab}(\lambda)}{\lambda^j}=\frac{a_j+b_j}{2}+O(|\lambda|)$$ and
letting $\lambda\to 0$ we obtain (1).
\end{proof}

As a consequence of Proposition 2 we obtain that any balanced
domain, which can be exhausted by domains biholomorphic to convex
domains is convex itself.

Note also that the condition $2k_{m+1}>k_n$ is essential as the
following simple example shows. The $(1,2)$-balanced domain
$$D=\{z\in\Bbb C^2:|z_1|^2+|z_2+z_1^2|<1\}$$ is not convex, but it
is biholomorphic to the $(1,2)$-balanced convex domain $$G
=\{z\in\Bbb C^2:|z_1|^2+|z_2|<1\}.$$

Now we are ready to prove our main result. To do this, we shall
apply Proposition 2 and the Cohn critertion which states (see e.g.
\cite{Rah-Sch}):

All the roots of a polynomial $\ds f(\zeta)=\sum_{j=0}^n
a_j\zeta^{n-j},$ $n\ge 2,$ $a_0\neq 0,$ belong to $\Bbb D$ if and
only if $|a_0|>|a_n|$ and all the roots of the polynomial
$$f^\star(\zeta)=\frac{\overline{a_0}f(\zeta)-a_n\overline
f(1/\overline\zeta)}{\zeta}$$ belong to $\Bbb D.$

\begin{proposition} The symmetrized n-disc $\Bbb G_n,$ $n\ge
3,$ cannot be exhausted by domains biholomorphic to convex
domains.
\end{proposition}

\begin{proof} Note that $\Bbb G_n$ is a $(1,2,\dots,n)$-balanced domain.
Hence, by Proposition 2, it is enough to show that if $\ds
m=\left[\frac{n}{2}\right]$, then the set $G_n$ of points
$(a_{m+1},\dots,a_n)\in\Bbb C^{n-m}$ such that all the zeros of
the polynomial $f_n(z)=z^n+\sum_{j=m+1}^n a_j\zeta^{n-j}$ belong
to $\Bbb D$ is not convex.

We shall first settle the cases $n=3$ and $n=4,$ and then we shall
reduce the general case to them.

{\it The case $n=3.$} For $f_3(\zeta)=\zeta^3+p\zeta+q$ one has
$$f_3^\star(\zeta)=\frac{f_3(\zeta)-q\overline
f_3(1/\overline\zeta)}{\zeta}=(1-|q|^2)\zeta^2-\overline
pq\zeta+p$$ and
$$f_3^{\star\star}(\zeta)=\frac{(1-|q|^2)f_3^\star(\zeta)-
p\overline{f_3^\star}(1/\overline\zeta)}{\zeta}=
((1-|q|^2)^2-|p|^2)\zeta-\overline pq(1-|q|^2)+p^2\overline q.$$
It follows from the Cohn criterion that $$G_3=\{(p,q)\in\Bbb
C^2:|q|<1,\ r(p,q)<0\},$$ where $$r(p,q)=|\overline
pq(1-|q|^2)-p^2\overline q|+|p|^2-(1-|q|^2)^2.$$ It is easy to see
that if $q'\in(-1,1)$ and $p'=1-q'^2,$ then $\ds
(p_1,q_1)=\left(p'e^{\frac{2\pi i}{3}},q'\right)$ and $\ds
(p_2,q_2)=\left(p'e^\frac{\pi i}{3},q'e^\frac{\pi i}{2},\right)$
are boundary points of $D,$ since $r(p',q')=0$ and $r(p',q)<0$ if
$p\in(|q'|-1,p').$ Then for
$$(p_0,q_0)=\left(\frac{p_1+q_1}{2},\frac{p_2+q_2}{2}\right)=
\left(p'\cos\frac{\pi}{6}e^\frac{\pi
i}{2},q'\cos\frac{\pi}{4}e^{\frac{\pi i}{4}}\right)$$ one has
$$|\overline{p_0}q_0(1-|q_0|^2)-p_0^2\overline{q_0}|=
|p_0q_0|(1-|q_0|^2+|p_0|).$$ Therefore
$$r(p_0,q_0)=(1-|q_0|^2+|p_0|)(1+|q_0|)(|p_0|+|q_0|-1).$$ So
$r(p_0,q_0)>0$ if and only if $|p_0|+|q_0|>1.$ For $\ds
q'=\frac{1}{2}$ it follows that
$$|p_0|+|q_0|=\frac{3\sqrt3+2\sqrt2}{8}>1.$$ Thus
$(p_0,q_0)\not\in\overline G_3$ and hence $G_3$ is not a convex
set.

{\it The case $n=4.$} Similar calculations as in the previous case
lead to $$G_4=\{(p,q)\in\Bbb C^2:|p|+|q|^2<1,\ s(p,q)<0\},$$ where
$$s(p,q)=(1-|q|^2)|\overline pq((1-|q|^2)^2-|p|^2)-p^3\overline
q^2|+|p|^4|q|^2-((1-|q|^2)^2-|p|^2)^2.$$ It is easy to see that if
$q'\in[0,1)$ and $p'=(1-q')\sqrt{1+q'},$ then $\ds
(p_1,q_1)=(p'e^\frac{\pi i}{2},q')\in\partial D$ and $\ds
(p_2,q_2)=(p'e^\frac{\pi i}{4},q'e^\frac{\pi i}{3})\in\partial D,$
since $s(p',q')=0$ and $s(p',q)<0$ if $p\in(-p',p').$ Then for
$$(p_0,q_0)=\left(\frac{p_1+q_1}{2},\frac{p_2+q_2}{2}\right)=
\left(p'\cos\frac{\pi}{8}e^\frac{3\pi
i}{8},q'\cos\frac{\pi}{6}e^\frac{\pi i}{6}\right)$$ one has
$$|\overline{p_0}q_0((1-|q_0|^2)^2-|p_0|^2)-p_0^3\overline
{q_0}^2|=|p_0q_0|((1-|q_0|^2)^2-|p_0|^2+|p_0|^2|q_0|).$$ Therefore
$$s(p_0,q_0)=(1-|q_0|^2)((1-|q_0|^2)(1+|q_0|)-|p_0|^2)(1+|p_0|-|q_0|^2)
(|p_0|+|q_0|-1).$$ So $s(p_0,q_0)>0$ if and only if
$|p_0|+|q_0|>1.$ For $\ds q'=\frac{2}{5}$ it follows that
$$|p_0|+|q_0|=\frac{1}{10}\left(3\sqrt{\frac{7(2+\sqrt2)}{5}}+2\sqrt
3\right)>1.$$ Thus $(p_0,q_0)\not\in\overline{G_4}$ and hence
$G_4$ is not a convex set.

{\it The case $n\ge 5.$} Let $j=\{0,1,2\}.$ Observe that the
non-convex set $G_3$ coincides with the set of points
$(p,q)\in\Bbb C^2$ such that all the zeros of the polynomial
$z^jf_3(z^k),$ $k\ge 1,$ belong to the unit disc. It follows that
if $n=3k+2$ and $k\ge 3$, $n=3k+1$ and $k\ge 2,$ or $n=3k$ and
$k\ge 1,$ then $G_3$ can be considered as an intersection of $G_n$
and a complex hyperplane. Therefore $G_n$ is not a convex set in
these cases.

In the remaining cases $n=5$ and $n=8$ it is enough to observe
that the non-convex set $G_4$ coincides with the set of points
$(p,q)\in\Bbb C^2$ such that all the zeros either of the
polynomials $\zeta f_4(\zeta)$ and $f_4(\zeta^2)$ belong to the
unit disc and then to complete the proof as above.\end{proof}

\end{document}